\documentclass{amsproc}
\usepackage{amsfonts}

\usepackage{graphicx}
\usepackage{amscd}

\theoremstyle{plain}

\newtheorem{corollary}{Corollary}

\newtheorem{definition}{Definition}

\newtheorem{lemma}{Lemma}

\newtheorem{proposition}{Proposition}
\newtheorem{remark}{Remark}

\newtheorem{theorem}{Theorem}
\numberwithin{equation}{section}

\input{tcilatex}

\begin{document}
\title[]{Gevrey vectors of multi-quasi-elliptic systems }
\author{Chikh BOUZAR}
\address{D\'{e}partement de Math\'{e}matiques. Universit\'{e} d'Oran Esenia.
Alg\'{e}rie.}
\email{bouzarchikh@hotmail.com}
\author{Rachid CHAILI}
\date{}

\begin{abstract}
We show that the multi--quasi-ellipticity is a necessary and sufficient
condition for the property of elliptic iterates to be hold for
multi-quasi-homogenous differential \ operators.
\end{abstract}

\maketitle

\section{Introduction}

The aim of this work is to prove the property of elliptic iterates for
multi-quasi-elliptic systems of differential operators in generalized Gevrey
spaces $G^{\mathcal{F},s}\left( \Omega \right) ,$ where $\mathcal{F}$
denotes the Newton's polyhedron of the system $\left( P_{j}\right)
_{j=1}^{N} $. The property of elliptic iterates for the system $\left(
P_{j}\right) _{j=1}^{N}$\ in the generalized Gevrey classes $G^{\mathcal{F}%
,s}\left( \Omega \right) $ means the inclusion 
\begin{equation*}
G^{s}\left( \Omega ,\left( P_{j}\right) _{j=1}^{N}\right) \subset G^{%
\mathcal{F},s}\left( \Omega \right)
\end{equation*}

Let $P_{j}\left( x,D\right) =\sum_{\alpha }a_{j\alpha }\left( x_{0}\right)
D^{\alpha },j=1,..,N,\;$denoted $\left( P_{j}\right) _{j=1}^{N},$ be linear
differential operators with $C^{\infty }$ coefficients in a open subset $%
\Omega $\ of $\mathbb{R}^{n}.$

\begin{definition}
Fix $x_{0}\in \Omega ,$ we define the Newton's polyhedron of the system $%
\left( P_{j}\right) _{j=1}^{N}$ at the point $x_{0}$, noted $\mathcal{F}%
\left( x_{0}\right) ,$\ as the convex hull of the set $\left\{ 0\right\}
\cup \left\{ \alpha \in \mathbb{N}^{n},\exists j\in \left\{ 1,..,N\right\}
;\;a_{j\alpha }\left( x_{0}\right) \neq 0\right\} .$ A Newton's polyhedron $%
\mathcal{F}$ is said regular if there exists a finite set $Q\left( \mathcal{F%
}\right) \subset \left( \mathbb{R}_{+}^{\ast }\right) ^{n}$\ such that 
\begin{equation*}
\mathcal{F}\mathbb{=}\underset{q\in Q\left( \mathcal{F}\right) }{\mathbb{%
\cap }}\left\{ q\in \mathbb{R}_{+}^{n},\;<\alpha ,q>\leq 1\right\}
\end{equation*}
\end{definition}

If $S\left( \mathcal{F}\right) $ designs the set of vertices of $\mathcal{F}$%
, we put

\begin{equation*}
V\left( \xi \right) =\sum_{\alpha \in S\left( \mathcal{F}\right) }\left| \xi
^{\alpha }\right| ,
\end{equation*}
and 
\begin{eqnarray*}
k\left( \alpha ,\mathcal{F}\right) &=&\inf \left\{ t>0,\;t^{-1}\alpha \in 
\mathcal{F}\right\} ,\text{ }\alpha \in \mathbb{R}_{+}^{n} \\
\mu _{j}\left( \mathcal{F}\right) &=&\max_{q\in Q\left( \mathcal{F}\right)
}q_{j}^{-1},\;j=1,..,n \\
\mu \left( \mathcal{F}\right) &=&\max_{1\leq j\leq n}\mu _{j}\left( \mathcal{%
F}\right) \\
\theta \left( \mathcal{F}\right) &=&\left( \frac{\mu \left( \mathcal{F}%
\right) }{\mu _{1}\left( \mathcal{F}\right) },..,\frac{\mu \left( \mathcal{F}%
\right) }{\mu _{n}\left( \mathcal{F}\right) }\right)
\end{eqnarray*}

\begin{definition}
Let $\mathcal{F}$\ be a regular polyhedron and $s\in \mathbb{R}_{\mathbb{+}%
},\;$we denote $G^{\mathcal{F},s}\left( \Omega \right) $ the space of $u\in
C^{\infty }\left( \Omega \right) $\ such that $\forall
H\;compact\;of\;\Omega ,\exists C>0,\forall \alpha \in \mathbb{N}^{n},$%
\begin{equation}
\sup_{H}\left| D^{\alpha }u\right| \leq C^{\left| \alpha \right| +1}\left[
\Gamma \left( \mu \left( \mathcal{F}\right) k\left( \alpha ,\mathcal{F}%
\right) +1\right) \right] ^{s}  \label{1.1}
\end{equation}
\end{definition}

\begin{remark}
We can take $\sup_{H}\left| D^{\alpha }u\right| $ or $\left\| D^{\alpha
}u\right\| _{L^{2}\left( H\right) }$ in the definition, according to Sobolev
imbedding theorems.
\end{remark}

\begin{definition}
The system $\left( P_{j}\right) _{j=1}^{N}$ is said multi-quasi-elliptic in $%
\Omega $\ if

1) $\mathcal{F}\left( x\right) =\mathcal{F}$\ does not depend of $x\in
\Omega $

2) $\mathcal{F}$ is regular

3) $\forall x\in \Omega ,\exists C>0,\exists R\geq 0,\forall \xi \in \mathbb{%
R}^{n},\left| \xi \right| \geq R,\;$%
\begin{equation*}
V\left( \xi \right) \leq C\sum\limits_{j=1}^{N}\left| P_{j}\left( x,\xi
\right) \right| \;
\end{equation*}
\end{definition}

\begin{definition}
Let $\left( P_{j}\right) _{j=1}^{N}$ be a system of linear differential
operators satisfying the conditions 1) and 2) of definition 3 and $s\in 
\mathbb{R}_{+},$ we denote $G^{s}\left( \Omega ,\left( P_{j}\right)
_{j=1}^{N}\right) $ the space of $u\in C^{\infty }\left( \Omega \right) $\
such that $\forall H\;compact\;of\;\Omega ,\exists C>0,\forall l\in \mathbb{N%
},1\leq i_{l}\leq N,$%
\begin{equation}
\left\| P_{i_{1}}...P_{i_{l}}u\right\| _{L^{2}\left( H\right) }\leq
C^{l+1}\left( l!\right) ^{s\mu \left( \mathcal{F}\right) }  \label{1.2}
\end{equation}
\end{definition}

The aim of this work is to show the following theorem.

\begin{theorem}
Let $\Omega $\ be an open subset of $\mathbb{R}^{n},\sigma >s\geq 1$ and $%
\left( P_{j}\left( x,D\right) \right) _{j=1}^{N}$ be a system of linear
differential operators with $G^{\theta \left( \mathcal{F}\right) ,\sigma
}\left( \Omega \right) $ coefficients, then 
\begin{equation*}
\left( P_{j}\right) _{j=1}^{N}\text{ is multi-quasi-elliptic in }\Omega
\Longleftrightarrow G^{s}\left( \Omega ,\left( P_{j}\right)
_{j=1}^{N}\right) \subset G^{\mathcal{F},s}\left( \Omega \right)
\end{equation*}
\end{theorem}

For differential operators with constant coefficients we have shown in \cite
{B-Ch2}\ a more general result.

\section{Sufficient condition}

Instead of $Q\left( \mathcal{F}\right) ,k\left( \mathcal{F},\alpha \right)
,\mu \left( \mathcal{F}\right) ,\mu _{j}\left( \mathcal{F}\right) ,\theta
\left( \mathcal{F}\right) $ we write, respectively, $Q,k\left( \alpha
\right) ,\mu ,\mu _{j},\theta .$ Denote $\mathcal{K}=\left\{ k=k\left(
\alpha \right) :\alpha \in \mathbb{N}^{n}\right\} $ and $\omega $ any open
subset of $\mathbb{R}^{n}.$ If $u\in C^{\infty }\left( \omega \right) $ and $%
k\in \mathcal{K},$ define 
\begin{equation*}
\left| u\right| _{k,\omega }=\sum\limits_{k\left( \alpha \right) =k}\left\|
D^{\alpha }u\right\| _{L^{2}\left( \omega \right) },
\end{equation*}
if $u\in C_{0}^{\infty }\left( \mathbb{R}^{n}\right) $ we write $\left|
u\right| _{k}.$ For $\rho >0$, we denote $B_{\rho }=\left\{ x\in \mathbb{R}%
^{n},\;\sum\limits_{j=1}^{n}x_{j}^{2\mu _{j}/\mu }<\rho ^{2}\right\} .$

Let $\left( P_{j}\right) _{j=1}^{N}$ be a system of linear differential
operators with coefficients defined in an open neighborhood $\Omega $ of the
origin satisfying the following conditions \newline
i) The system $\left( P_{j}\right) _{j=1}^{N}$ is multi-quasi-elliptic in $%
\Omega .$\newline
ii) The coefficients $a_{j\alpha }\in G^{\theta ,s}\left( \Omega \right)
,\forall \alpha \in \mathbb{N}^{n},\forall j\in \left\{ 1,..,N\right\} .$

We define for $j=1,...,N,$ and $h\in \mathbb{N},$ 
\begin{equation*}
P_{j}^{h}\left( x,D\right) =\underset{h\;times}{\underbrace{P_{j}\left(
x,D\right) \circ \cdot \cdot \circ P_{j}\left( x,D\right) }}.
\end{equation*}

From the multi-quasi-ellipticity of the system $\left( P_{j}\right)
_{j=1}^{N}$ and following the proof of the lemma 3.4 of \cite{Z}, we easily
obtain

\begin{lemma}
There exists $\rho _{0}>0,$ $\forall \varepsilon \in ]0,\frac{1}{v\left(
n\right) }[,$($v\left( n\right) $ denote the number of elements of $\mathcal{%
K}\cap \lbrack 0,n[$), we have\newline
$\exists C_{1}>0,\exists C_{2}\left( \varepsilon \right) >0,\forall \delta
\in ]0,1[,\forall \rho >0,B_{\rho +\delta }\subset B_{\rho _{0}},\forall
u\in C^{\infty }\left( B_{\rho _{0}}\right) ,\forall p\geq n,$%
\begin{eqnarray}
\left| u\right| _{p+1,B_{\rho }} &\leq &C_{1}\left(
\sum\limits_{j=1}^{N}\left| P_{j}^{n}(x,D)u\right| _{p-n+1,B_{\rho +\delta
}}+\varepsilon \left| u\right| _{p+1,B_{\rho +\delta }}+\left( \varepsilon
\delta \right) ^{-n\mu }\left| u\right| _{p-n+1,B_{\rho +\delta }}\right. 
\notag \\
&&\left. +\sum\limits_{h=0}^{p}\left( \frac{\left( p+1\right) !}{h!}\right)
^{s\mu }C_{2}\left( \varepsilon \right) ^{p+1-h}\left| u\right| _{h,B_{\rho
+\delta }}\right) \;,  \label{2.1}
\end{eqnarray}
and for $p\leq n$, we have 
\begin{equation}
\left| u\right| _{p+1,B_{\rho }}\leq C_{1}\left( \sum\limits_{j=1}^{N}\left|
P_{j}^{n}(x,D)u\right| _{p-n+1,B_{\rho +\delta }}+\varepsilon \left|
u\right| _{p+1,B_{\rho +\delta }}+\left( \varepsilon \delta \right)
^{-\left( p+1\right) \mu }\left| u\right| _{0,B_{\rho +\delta }}\right)
\label{2.2}
\end{equation}
\end{lemma}

Let $\lambda >0$ and $R>0,$ for $p\in \mathbb{N},\;$we set 
\begin{equation*}
\sigma _{p}\left( u,\lambda \right) =\sigma _{p}\left( u,\lambda ,R\right)
=\left( p!\right) ^{-s\mu }\lambda ^{-p}\sup_{R/2\leq \rho <R}\left( R-\rho
\right) ^{p\mu }\left| u\right| _{p,B_{\rho }}
\end{equation*}

\begin{lemma}
Let $\rho _{0}$ be as in the lemma1 and let $0<R<1$ such that $\overline{B}%
_{R}\subset B_{\rho _{0}},$ then $\exists \lambda _{0}>0\left( \lambda _{0}%
\text{ depends only of }R\text{ and }\left( P_{j}\right) _{j=1}^{N}\right)
,\forall u\in C^{\infty }\left( B_{\rho _{0}}\right) ,\forall \lambda \geq
\lambda _{0},\forall p\geq n,$%
\begin{equation}
\sigma _{p+1}\left( u,\lambda \right) \leq \left[ \left( p-n+2\right)
...\left( p+1\right) \right] ^{-s\mu }\sum\limits_{j=1}^{N}\sigma
_{p-n+1}\left( P_{j}^{n}u,\lambda \right) +\sum\limits_{h=0}^{p}\sigma
_{h}\left( u,\lambda \right) \;,  \label{2.3}
\end{equation}
and for $p\leq n-1,$%
\begin{equation}
\sigma _{p+1}\left( u,\lambda \right) \leq \left( p+1\right) !^{-s\mu
}\sum\limits_{j=1}^{N}\sigma _{0}\left( P_{j}^{p+1}u,\lambda \right) +\sigma
_{0}\left( u,\lambda \right)  \label{2.4}
\end{equation}
\end{lemma}

\begin{proof}
Let $p\geq n,$ multiplying both sides of inequality $\left( 2.1\right) $ by $%
\left( p+1\right) !^{-s\mu }\lambda ^{-p-1}\left( R-\rho \right) ^{p\mu },$
taking $\delta =\frac{R-\rho }{p-n+2}$ and then passing to the $\sup $ over $%
\rho \in \lbrack R/2,R[$, we obtain 
\begin{equation*}
\sigma _{p+1}\left( u,\lambda \right) \leq C_{1}\left( I_{1}+\varepsilon
I_{2}+\varepsilon ^{-n\mu }I_{3}+I_{4}\right) \;,
\end{equation*}
where 
\begin{eqnarray*}
I_{1} &=&\sum\limits_{j=1}^{N}\frac{\left( R-\rho \right) ^{\left(
p+1\right) \mu }}{\left( p+1\right) !^{s\mu }\lambda ^{p+1}}\left|
P_{j}^{n}u\right| _{p-n+1,B_{\rho +\delta }} \\
&\leq &\sum\limits_{j=1}^{N}\left( \frac{\left( p-n+1\right) !}{\left(
p+1\right) !}\right) ^{s\mu }\frac{e^{\mu }}{\lambda ^{n}}\sigma
_{p-n+1}\left( P_{j}^{n}u,\lambda \right) \\
I_{2} &=&\frac{\left( R-\rho \right) ^{\left( p+1\right) \mu }}{\left(
p+1\right) !^{s\mu }\lambda ^{p+1}}\left| u\right| _{p+1,B_{\rho +\delta
}}\leq \left( 2^{n}e\right) ^{\mu }\sigma _{p+1}\left( u,\lambda \right) \\
I_{3} &=&\frac{\left( R-\rho \right) ^{\left( p+1\right) \mu }}{\left(
p+1\right) !^{s\mu }\lambda ^{p+1}}\delta ^{-n\mu }\left| u\right|
_{p-n+1,B_{\rho +\delta }}\leq \frac{e^{\mu }}{\lambda ^{n}}\sigma
_{p-n+1}\left( u,\lambda \right) \\
I_{4} &=&\frac{\left( R-\rho \right) ^{\left( p+1\right) \mu }}{\left(
p+1\right) !^{s\mu }\lambda ^{p+1}}\sum\limits_{h=0}^{p}\left( \frac{\left(
p+1\right) !}{h!}\right) ^{s\mu }C_{2}\left( \varepsilon \right)
^{p+1-h}\left| u\right| _{h,B_{\rho +\delta }} \\
&\leq &\frac{e^{\mu }C_{2}\left( \varepsilon \right) }{\lambda }%
\sum\limits_{h=0}^{p}\left( \frac{C_{2}\left( \varepsilon \right) }{\lambda }%
\right) ^{p-h}\sigma _{h}\left( u,\lambda \right)
\end{eqnarray*}
By a suitable choice of $\varepsilon $, we find 
\begin{eqnarray*}
\sigma _{p+1}\left( u,\lambda \right) &\leq &\left( \frac{\left(
p-n+1\right) !}{\left( p+1\right) !}\right) ^{s\mu }\frac{\widetilde{C}_{1}}{%
\lambda ^{n}}\sum\limits_{j=1}^{N}\sigma _{p-n+1}\left( P_{j}^{n}u,\lambda
\right) +\frac{\widetilde{C}_{2}}{\lambda ^{n}}\sigma _{p-n+1}\left(
u,\lambda \right) \\
&&+\frac{\widetilde{C}_{3}}{\lambda }\sum\limits_{h=0}^{p}\left( \frac{%
\widetilde{C}_{4}}{\lambda }\right) ^{p-h}\sigma _{h}\left( u,\lambda \right)
\end{eqnarray*}
It suffices to take $\lambda _{0}=$ $\widetilde{C}_{1}+\widetilde{C}_{2}+%
\widetilde{C}_{3}+\widetilde{C}_{4}$ to get $\left( 2.3\right) .$ For the
inequality $\left( 2.4\right) $ we multiply both sides of inequality $\left(
2.2\right) $ by $\frac{\left( R-\rho \right) ^{\left( p+1\right) \mu }}{%
\left( p+1\right) !^{s\mu }\lambda ^{p+1}}$, take $\delta =\frac{R-\rho }{2}$
and then we follow the same procedure as for $\left( 2.3\right) $
\end{proof}

\begin{lemma}
Let $\rho _{0},R$ and $\lambda _{0}$ be as in lemma 2, then $\forall u\in
C^{\infty }\left( B_{\rho _{0}}\right) ,\;\forall \lambda \geq \lambda
_{0},\;\forall p\in \mathbb{N},$%
\begin{equation}
\sigma _{p+1}\left( u,\lambda \right) \leq 2^{p+1}\sigma _{0}\left(
u,\lambda \right) +\sum\limits_{l=1}^{p+1}2^{p+1-l}C_{p+1}^{l}\frac{1}{%
\left( l!\right) ^{s\mu }}\sum\limits_{\substack{ 1\leq j\leq l  \\ 1\leq
i_{j}\leq N}}\sigma _{0}\left( P_{i_{1}}..P_{i_{l}}u,\lambda \right)
\label{2.5}
\end{equation}
\end{lemma}

\begin{proof}
It's obtained by recurrence over $p$ as in lemma 3.2 of \cite{B-Ch1}.
\end{proof}

Our first result is the following theorem.

\begin{theorem}
Let $\Omega $\ be an open subset of $\mathbb{R}^{n},\,s\geq 1$ and let $%
\left( P_{j}(x,D)\right) _{j=1}^{N}$ be a system of linear differential
operators with $G^{\theta ,s}\left( \Omega \right) $ coefficients$,$ then 
\begin{equation*}
\left( P_{j}\right) _{j=1}^{N}\text{ is multi-quasi-elliptic in }\Omega
\Rightarrow G^{s}\left( \Omega ,\left( P_{j}\right) _{j=1}^{N}\right)
\subset G^{\mathcal{F},s}\left( \Omega \right)
\end{equation*}
\end{theorem}

\begin{proof}
It is sufficient to check $\left( 1.1\right) $ in a neighborhood of every
point $x$ of $\Omega .$ Let $x$ be the origin, then there exist $\rho
_{0},\lambda _{0}$ and $R$ such that the lemmas hold.\newline
Let $u\in G^{s}\left( \Omega ,\left( P_{j}\right) _{j=1}^{N}\right) $, then
there exists $C_{1}>0$ such that 
\begin{equation*}
\sigma _{0}\left( P_{i_{1}}..P_{i_{l}}u,\lambda _{0}\right) \leq
C_{1}^{l+1}\left( l!\right) ^{s\mu },\;\forall l\in \mathbb{N},
\end{equation*}
hence from $\left( 2.5\right) ,$ we obtain 
\begin{equation*}
\sigma _{p+1}\left( u,\lambda _{0}\right) \leq C_{1}\left( 2+NC_{1}\right)
^{p+1},\;\forall p\in \mathbb{N},
\end{equation*}
which gives 
\begin{equation}
\left| u\right| _{p+1,B_{R/2}}\leq \left( p+1\right) !^{s\mu }C_{2}^{\left(
p+1\right) \mu +1}\;,\forall p\in \mathbb{N}  \label{2.6}
\end{equation}
Let $k\in \mathcal{K},$ if $k\geq n-1,$ then $k=p+r,$ where $p=\left[ k-n%
\right] +1$ and $r\in \lbrack n-1,n[.$ The interpolation inequality from 
\cite[lemme 2.3]{Z} gives 
\begin{equation*}
\left| u\right| _{k,B_{R/2}}\leq \varepsilon \left( \left| u\right|
_{p+n,B_{R/2}}+\left| u\right| _{p+n-1,B_{R/2}}+\sum\limits_{\substack{ t\in 
\mathcal{K}  \\ p\leq t<p+n-1}}\left| u\right| _{t,B_{R/2}}\right)
+C\varepsilon ^{-\frac{r}{n-r}}\left| u\right| _{p,B_{R/2}}
\end{equation*}
If we set $\varepsilon =\left( \frac{\Gamma \left( k+1\right) }{\Gamma
\left( p+n+1\right) }\right) ^{s\mu },$ then we obtain, with $\left(
2.6\right) ,$%
\begin{equation*}
\left| u\right| _{k,B_{R/2}}\leq C_{2}^{k\mu +1}\left( \Gamma \left(
k+1\right) \right) ^{s\mu }\left( 2+CC_{3}\right) +\sum\limits_{\substack{ %
t\in \mathcal{K}  \\ p\leq t<p+n-1}}\left| u\right| _{t,B_{R/2}},
\end{equation*}
where $C_{3}$ is the constant of the inequality 
\begin{equation*}
\left( \frac{\Gamma \left( p+n+1\right) }{\Gamma \left( p+r+1\right) }%
\right) ^{\frac{r}{n-r}}\leq C_{3}^{1/\mu }\frac{\Gamma \left( p+r+1\right) 
}{\Gamma \left( p+1\right) },\;p\in \mathbb{N},\;r\in \mathcal{K}\cap
\lbrack 0,n[
\end{equation*}
Applying again the above interpolation inequality to $\sum\limits_{\substack{
t\in \mathcal{K}  \\ h\leq t<h+n-1}}\left| u\right| _{t,B_{R/2}},$ $%
h=0,..,p, $ we find 
\begin{equation*}
\left| u\right| _{k,B_{R/2}}\leq C_{2}^{k\mu +1}\left( \Gamma \left(
k+1\right) \right) ^{s\mu }\left( 2+CC_{3}\right)
\sum\limits_{j=0}^{p}v\left( n\right) ^{j},
\end{equation*}
where $v\left( n\right) $ is the number of elements of $\mathcal{K}\cap
\lbrack 0,n[.$\newline
Similarly, if $k\leq n-1,$ we find 
\begin{equation*}
\left| u\right| _{k,B_{R/2}}\leq C_{2}^{k\mu +1}\left( \Gamma \left(
k+1\right) \right) ^{s\mu }\left( 1+CC_{3}\right)
\end{equation*}
Hence for all $k\in \mathcal{K}$, we have 
\begin{equation}
\left| u\right| _{k,B_{R/2}}\leq \widetilde{C}^{k\mu +1}\left( \Gamma \left(
k+1\right) \right) ^{s\mu }  \label{2.7}
\end{equation}
By the imbedding theorem of anisitropic Sobolev spaces, if $l\in \mathcal{K}$
and $l>\frac{2}{k\left( e\right) },$ where $k\left( e\right) =\underset{q\in
Q}{\max }\sum\limits_{j=1}^{n}q_{j},$ then $\exists C>0,\forall \alpha \in 
\mathbb{N}^{n},$%
\begin{equation*}
\sup_{B_{R/2}}\left| D^{\alpha }u\right| \leq C\left| D^{\alpha }u\right|
_{l,B_{R/2}},
\end{equation*}
consequently from $\left( 2.7\right) $ we get the estimate $\left(
1.1\right) $
\end{proof}

As a consequence of this theorem, we obtain a result of Gevrey
hypo-ellipticity for multi-quasi-elliptic systems.

\begin{corollary}
Under the assumptions of theorem 2, the following propositions are equivalent

i) $u\in \mathcal{D}^{\prime }\left( \Omega \right) ,$ $P_{j}u\in G^{%
\mathcal{F},s}\left( \Omega \right) ,\;\forall j=1,..,N.$

ii) $u\in G^{\mathcal{F},s}\left( \Omega \right) .$
\end{corollary}

\section{Necessary condition}

In this section we prove the reciprocal of the theorem 2. For this aim we
need a characterization of the multi-quasi-ellipticity for the system $%
\left( P_{j}(x,D)\right) _{j=1}^{N}$, known in the case of a scalar operator
in \cite{G-V}.

\begin{proposition}
A system $\left( P_{j}\right) _{j=1}^{N},$ satisfying 1) and 2) of
Definition 2, is multi-quasi-elliptic in $\Omega $\ if and only if for any $%
x\in \Omega ,$ $\forall q\in Q,$ 
\begin{equation*}
\sum\limits_{j=1}^{N}\left| P_{jq}\left( x,\xi \right) \right| \neq
0,\;\;\forall \xi \in \mathbb{R}^{n},\xi _{1}...\xi _{n}\neq 0\;,
\end{equation*}
where $P_{jq}$ is the $q$-quasi-homogenous part of $P_{j},$ i.e. 
\begin{equation*}
P_{jq}\left( x,\xi \right) =\sum\limits_{<\alpha ,q>=1}a_{j\alpha }\left(
x\right) \xi ^{\alpha }
\end{equation*}
\end{proposition}

\begin{theorem}
Let\textit{\ }$\Omega $\textit{\ be an open subset of }$\mathbb{R}^{n}$%
\textit{\ and }$P_{j}\left( x,D\right) ,$\textit{\ }$j=1,..,N,$\textit{\ be
differential operators with }$G^{\theta ,\sigma }\left( \Omega \right) $%
\textit{\ coefficients, if }$s>\sigma \geq 1,$ then\textit{\ } 
\begin{equation*}
G^{s}\left( \Omega ,\left( P_{j}\right) _{j=1}^{N}\right) \subset G^{%
\mathcal{F},s}\left( \Omega \right) \Rightarrow \left( P_{j}\right)
_{j=1}^{N}\text{\textit{\ }is multi-\textit{quasi-elliptic in} }\Omega
\end{equation*}
\end{theorem}

\begin{proof}
Assume that the system $\left( P_{j}\right) _{j=1}^{N}$ is not
multi-quasi-elliptic, then there exists $x_{0}\in \Omega ,\exists q\in Q$
and $\exists \xi _{0}\in S^{n-1},\xi _{0,1}...\xi _{0,n}\neq 0,$ such that 
\begin{equation}
P_{jq}\left( x_{0},\xi _{0}\right) =0,\quad \forall j=1,..,N  \label{3.1}
\end{equation}

We shall construct a function $u\in G^{s}\left( \Omega ,\left( P_{j}\right)
_{j=1}^{N}\right) $ and $u\notin G^{\mathcal{F},s}\left( \Omega \right) $,
wich contradicts the hypothesis. Choose $\varepsilon $ such that 
\begin{equation*}
0<\varepsilon \leq \frac{\mu \left( s-\sigma \right) }{2\mu s-\sigma }<\frac{%
1}{2}\text{ and }\varepsilon \leq \underset{<\beta ,q><1}{\min }\mu \left(
1-<\beta ,q>\right) \quad ,
\end{equation*}
and put $\eta =\dfrac{1-\varepsilon /\mu }{\mu s}.$\newline

let $\delta >0$ such that the ball $B_{0}=B\left( x_{0},2\delta \right) $ be
relatively compact in$\ \Omega ,$ and let $\varphi \in G^{q,\sigma \mu
}\left( \mathbb{R}^{n}\right) $ with compact support in $B\left( 0,2\delta
\right) $ and $\varphi \left( x\right) \equiv 1$ in $B\left( 0,\delta
\right) .$ The desired function is 
\begin{equation*}
u\left( x\right) =\int_{1}^{+\infty }\varphi \left[ r^{\varepsilon q}\left(
x-x_{0}\right) \right] e^{-r^{\eta }}e^{i<x-x_{0},r^{q}\xi _{0}>}dr\quad ,
\end{equation*}
where $r^{q}x=\left( r^{q_{1}}x_{1},r^{q_{2}}x_{2},..,r^{q_{n}}x_{n}\right)
. $

As in \cite{M}, we easly show, for sufficiently large $m,$ that 
\begin{equation*}
\left| D^{m\alpha }u\left( x_{0}\right) \right| >\frac{1}{2\eta }\left| \xi
_{0}^{m\alpha }\right| \Gamma \left( \frac{<m\alpha ,q>+1}{\eta }\right)
\end{equation*}
If we choose $\alpha $ such that $k\left( \alpha \right) =<\alpha ,q>,$ we
obtain 
\begin{equation*}
\left| D^{m\alpha }u\left( x_{0}\right) \right| >\frac{c\left( \mu \right)
^{<m\alpha ,q>}}{2\eta }\left| \xi _{0}^{m\alpha }\right| \left[ \Gamma
\left( \mu k\left( m\alpha \right) +1\right) \right] ^{\frac{1}{\mu \eta }%
}\;,
\end{equation*}
and, since $\frac{1}{\mu \eta }>s,$ then $u\notin G^{\mathcal{F},s}\left(
U_{0}\right) $ for any neighborhood $U_{0}$ of $x_{0}.$

Let show that $u\in G^{s}\left( \Omega ,\left( P_{j}\right)
_{j=1}^{N}\right) .$ Since the coefficients of $P_{j}$ are in $G^{\theta
,\sigma }\left( \Omega \right) \subset G^{q,\sigma \mu }\left( \Omega
\right) ,$ thus $\exists M>0,$ $\forall \alpha \in \mathbb{Z}_{+}^{n},$ $%
\forall \beta \in \mathbb{Z}_{+}^{n},$ $\forall x\in B_{0},$ $\forall r\geq
1,$ $\forall j=1,..,N,$%
\begin{equation}
\left| \left( D_{x}^{\beta }P_{j}^{(\alpha )}\right) \left( x,r^{q}\xi
_{0}\right) \right| \leq M^{\left| \beta \right| +1}\left[ \Gamma (<\beta
,q>+1)\right] ^{\sigma \mu }r^{1-<\alpha ,q>}  \label{3.2}
\end{equation}
In the other hand in view of $\left( 3.1\right) $ we have$\newline
\forall \delta >0,$ $\exists C_{1}>0,$ $\forall r\geq 1,$ $\forall x\in
\Omega ,\left| x-x_{0}\right| <2\delta r^{-\varepsilon /\mu },$ $\forall
j=1,..,N,$%
\begin{equation}
\left| P_{j}\left( x,r^{q}\xi _{0}\right) \right| \leq C_{1}r^{1-\varepsilon
/\mu }\quad  \label{3.3}
\end{equation}
Indeed 
\begin{eqnarray*}
\left| P_{j}\left( x,r^{q}\xi _{0}\right) \right| &\leq &\left| P_{jq}\left(
x,r^{q}\xi _{0}\right) \right| +\sum\limits_{<\beta ,q><1}r^{<\beta
,q>}\left| a_{j\beta }\left( x\right) \right| \left| \xi _{0}^{\beta }\right|
\\
&\leq &r\left| P_{jq}\left( x,\xi _{0}\right) -P_{jq}\left( x_{0},\xi
_{0}\right) \right| +\sum\limits_{<\beta ,q><1}r^{<\beta ,q>}\left|
a_{j\beta }\left( x\right) \right| \\
&\leq &r\sum\limits_{<\beta ,q>=1}\left| a_{j\beta }\left( x\right)
-a_{j\beta }\left( x_{0}\right) \right| +\sum\limits_{<\beta ,q><1}r^{<\beta
,q>}\left| a_{j\beta }\left( x\right) \right| \quad ,
\end{eqnarray*}
by the mean-value theorem and the choice of $\varepsilon ,$ we obtain 
\begin{equation*}
\left| P_{j}\left( x,r^{q}\xi _{0}\right) \right| \leq r\sum\limits_{<\beta
,q>=1}C2\delta r^{-\varepsilon /\mu }+\sum\limits_{<\beta ,q><1}C^{\prime
}r^{1-\varepsilon /\mu }\leq C_{1}r^{1-\varepsilon /\mu }
\end{equation*}
Since $\varphi \in G_{0}^{q,\sigma \mu }\left( \mathbb{R}^{n}\right) ,$ so $%
\exists C_{0}>0,$ $\exists L_{0}>0,$ $\forall \beta \in \mathbb{Z}_{+}^{n},$ 
$\forall x\in \mathbb{R}^{n}$%
\begin{equation}
\left| D^{\beta }\varphi \left( x\right) \right| \leq C_{0}L_{0}^{\left|
\beta \right| }\left[ \Gamma (<\beta ,q>+1)\right] ^{\sigma \mu }
\label{3.4}
\end{equation}
We choose $L_{0}\geq 2MC_{2}^{\sigma \mu },$ where $C_{2}$ is the constant
of the following inequality $\exists C>0,\;\forall \gamma \in \mathbb{Z}%
_{+}^{n},\forall \beta \in \mathbb{Z}_{+}^{n},$ 
\begin{equation}
\frac{\gamma !}{\beta !\left( \gamma -\beta \right) !}\leq C^{<\gamma -\beta
,q>}\frac{\Gamma (<\gamma ,q>+1)}{\Gamma (<\beta ,q>+1)\Gamma (<\gamma
-\beta ,q>+1)}  \label{3.5}
\end{equation}
In the sequel we will use the following properties of the gamma function. 
\begin{equation}
\Gamma (a+p)=\Gamma \left( a\right) a\left( a+1\right) \cdot \cdot
(a+p-1),\quad a>0,\text{ }p\in \mathbb{Z}_{+}^{\ast }  \label{3.6}
\end{equation}
$\forall \lambda >0,$ $\forall \tau >0,$ $\forall a\geq \omega ,$ $\forall
b\geq \omega $ $,\forall c\geq \omega ,\forall \sigma \geq 1$%
\begin{equation}
\lambda ^{a}\Gamma (b+c+1)^{\sigma }\tau ^{c}\leq 2^{\frac{\sigma }{\omega }%
}[\lambda ^{a+c}\Gamma (b+1)^{\sigma }+\Gamma (a+b+c+1)^{\sigma }\tau
^{a+c}].  \label{3.7}
\end{equation}

We need a convenient form of $P_{i_{k}}..P_{i_{1}}u$ for any integer $k\geq
0.$ By the generalized Leibniz formula $P_{j}\left( x,D\right) \left(
uv\right) =\sum\limits_{\alpha }\frac{1}{\alpha !}P_{j}^{\left( \alpha
\right) }uD^{\alpha }v\quad ,$ we obtain 
\begin{eqnarray*}
P_{j}u\left( x\right) &=&\int_{1}^{+\infty }\sum\limits_{<\alpha ,q>\leq 1}%
\frac{1}{\alpha !}D_{x}^{\alpha }\left( \varphi \left[ r^{\varepsilon
q}\left( x-x_{0}\right) \right] e^{-r^{\eta }}\right) P_{j}^{(\alpha
)}\left( e^{i<x-x_{0},r^{q}\xi _{0}>}\right) dr \\
&=&\int_{1}^{+\infty }\sum\limits_{<\alpha ,q>\leq 1}\frac{1}{\alpha !}%
P_{j}^{(\alpha )}\left( x,r^{q}\xi _{0}\right) D_{x}^{\alpha }\left( \varphi %
\left[ r^{\varepsilon q}\left( x-x_{0}\right) \right] \right) e^{-r^{\eta
}}e^{i<x-x_{0},r^{q}\xi _{0}>}dr \\
&=&\int_{1}^{+\infty }A_{j}\left( x,r\right) e^{-r^{\eta
}}e^{i<x-x_{0},r^{q}\xi _{0}>}dr\quad ,
\end{eqnarray*}
where 
\begin{equation*}
A_{j}\left( x,r\right) =\sum\limits_{<\alpha ,q>\leq 1}\frac{1}{\alpha !}%
P_{j}^{(\alpha )}\left( x,r^{q}\xi _{0}\right) D_{x}^{\alpha }\left( \varphi %
\left[ r^{\varepsilon q}\left( x-x_{0}\right) \right] \right)
\end{equation*}
For any integer $k\geq 0,$ $l\leq k$ and $1\leq i_{l}\leq N,$ we write 
\begin{equation*}
P_{i_{k}}..P_{i_{0}}u\left( x\right) =\int_{1}^{+\infty
}A_{i_{k}...i_{0}}\left( x,r\right) e^{-r^{\eta }}e^{i<x-x_{0},r^{q}\xi
_{0}>}dr\quad ,
\end{equation*}
where $P_{i_{0}}$ designs the identity operator, and 
\begin{equation}
\left\{ 
\begin{array}{c}
A_{i_{0}}\left( x,r\right) =\varphi \left[ r^{\varepsilon q}\left(
x-x_{0}\right) \right] \qquad \qquad \qquad \qquad \qquad \qquad \quad \\ 
A_{i_{k+1},i_{k}...i_{0}}\left( x,r\right) =\dsum\limits_{<\alpha ,q>\leq 1}%
\frac{1}{\alpha !}P_{i_{k+1}}^{(\alpha )}\left( x,r^{q}\xi _{0}\right)
D_{x}^{\alpha }A_{i_{k}...i_{0}}\left( x,r\right)
\end{array}
\right.  \label{3.8}
\end{equation}
To complete the proof we need the following

\begin{lemma}
$\exists L>0,$ $\forall k\in \mathbb{Z}_{+},$ $\forall \gamma \in \mathbb{Z}%
_{+}^{n},$ $\forall x\in B_{0},$ $\forall r\geq 1,$%
\begin{eqnarray}
\left| D_{x}^{\gamma }A_{i_{k}...i_{0}}\left( x,r\right) \right| &\leq
&C_{0}\left( L_{0}r^{\varepsilon }\right) ^{\left| \gamma \right|
}L^{k}\left( r^{\left( 1-\varepsilon /\mu \right) k}\left[ \Gamma (<\gamma
,q>+1)\right] ^{\sigma \mu }\right.  \notag \\
&&\left. +\left[ \Gamma (<\gamma ,q>+k+1)\right] ^{\sigma \mu
}r^{k\varepsilon \left( 2-1/\mu \right) }\right)  \label{3.9}
\end{eqnarray}
\end{lemma}

\textbf{Proof.} It's obtained by recurrence over $k.$ In fact for $k=0,$
this is true since $\varphi \in G_{0}^{q,\sigma \mu }\left( \mathbb{R}%
^{n}\right) ,$ i.e. we find $\left( 3.4\right) .$ So suppose that the
estimate $\left( 3.9\right) $ holds up to the order $k$ and check that it
remains valid at the order $k+1.$ Set\newline
$\lambda =r^{1-\varepsilon /\mu },$\newline
$\tau =r^{\varepsilon \left( 2-1/\mu \right) },\newline
S\left( k,\beta \right) =\lambda ^{k}\left[ \Gamma (<\beta ,q>+1)\right]
^{\sigma \mu }+\left[ \Gamma (<\beta ,q>+k+1)\right] ^{\sigma \mu }\tau
^{k}. $\newline
Then the estimate $\left( 3.9\right) $ is written as 
\begin{equation*}
\left| D_{x}^{\gamma }A_{i_{k}...i_{0}}\left( x,r\right) \right| \leq
C_{0}\left( L_{0}r^{\varepsilon }\right) ^{\left| \gamma \right|
}L^{k}S\left( k,\gamma \right)
\end{equation*}
Let $\omega =\underset{1\leq j\leq n}{\min }q_{j},$ from inequality $\left(
3.7\right) ,$ we have 
\begin{equation}
r^{1-\varepsilon /\mu }S\left( k,\beta \right) \leq 2^{\frac{\sigma \mu }{%
\omega }+1}S\left( k+1,\beta \right) \;,  \label{3.10}
\end{equation}
indeed 
\begin{eqnarray*}
r^{1-\varepsilon /\mu }S\left( k,\beta \right) &=&\lambda \left( \lambda
^{k} \left[ \Gamma (<\beta ,q>+1)\right] ^{\sigma \mu }+\left[ \Gamma
(<\beta ,q>+k+1)\right] ^{\sigma \mu }\tau ^{k}\right) \\
\ &=&\lambda ^{k+1}\left[ \Gamma (<\beta ,q>+1)\right] ^{\sigma \mu
}+\lambda \left[ \Gamma (<\beta ,q>+k+1)\right] ^{\sigma \mu }\tau ^{k} \\
\ &\leq &2^{\frac{\sigma \mu }{\omega }}\left( \lambda ^{k+1}\left[ \Gamma
(<\beta ,q>+1)\right] ^{\sigma \mu }+\left[ \Gamma (<\beta ,q>+\left(
k+1\right) +1)\right] ^{\sigma \mu }\tau ^{k+1}\right. \\
&&\ \ \ \ \ \ \left. +\lambda ^{k+1}\left[ \Gamma (<\beta ,q>+1)\right]
^{\sigma \mu }+\left[ \Gamma (<\beta ,q>+\left( k+1\right) +1)\right]
^{\sigma \mu }\tau ^{k+1}\right) \\
&\leq &2^{\frac{\sigma \mu }{\omega }+1}S\left( k+1,\beta \right)
\end{eqnarray*}
Similarly we find 
\begin{equation}
\lambda ^{1-<\alpha ,q>}\tau ^{<\alpha ,q>}S\left( k,\beta +\alpha \right)
\leq 2^{\frac{\sigma \mu }{\omega }+1}S\left( k+1,\beta \right) ,\quad
<\alpha ,q>\leq 1  \label{3.11}
\end{equation}
From $\left( 3.8\right) ,$ we have 
\begin{eqnarray*}
\left| D_{x}^{\gamma }A_{i_{k+1}...i_{0}}\left( x,r\right) \right| &=&\left|
\sum\limits_{<\alpha ,q>\leq 1}\frac{1}{\alpha !}D_{x}^{\gamma }\left(
P_{i_{k+1}}^{(\alpha )}\left( x,r^{q}\xi _{0}\right) D_{x}^{\alpha
}A_{i_{k}...i_{0}}\left( x,r\right) \right) \right| \\
\ &\leq &\left| \sum\limits_{<\alpha ,q>\leq 1}\sum\limits_{\beta \leq
\gamma }\frac{1}{\alpha !}\left( _{\beta }^{\gamma }\right) D_{x}^{\gamma
-\beta }P_{i_{k+1}}^{(\alpha )}\left( x,r^{q}\xi _{0}\right) D_{x}^{\alpha
+\beta }A_{i_{k}...i_{0}}\left( x,r\right) \right| \\
\ &\leq &I_{1}+I_{2}+I_{3}\quad ,
\end{eqnarray*}
where 
\begin{eqnarray*}
I_{1} &=&\left| P_{i_{k+1}}\left( x,r^{q}\xi _{0}\right) \right| \left|
D_{x}^{\gamma }A_{i_{k}...i_{0}}\left( x,r\right) \right| \\
I_{2} &=&\sum\limits_{\beta <\gamma }\left( _{\beta }^{\gamma }\right)
\left| D_{x}^{\gamma -\beta }P_{i_{k+1}}\left( x,r^{q}\xi _{0}\right)
\right| \left| D_{x}^{\beta }A_{i_{k}...i_{0}}\left( x,r\right) \right| \\
I_{3} &=&\sum\limits_{0<<\alpha ,q>\leq 1}\sum\limits_{\beta \leq \gamma }%
\frac{1}{\alpha !}\left( _{\beta }^{\gamma }\right) \left| D_{x}^{\gamma
-\beta }P_{i_{k+1}}^{(\alpha )}\left( x,r^{q}\xi _{0}\right) \right| \left|
D_{x}^{\alpha +\beta }A_{i_{k}...i_{0}}\left( x,r\right) \right|
\end{eqnarray*}
Since the $A_{i_{k}...i_{0}}$ are functions of compact supports in $B\left(
x_{0},2\delta r^{-\varepsilon /\mu }\right) ,$ so due to $\left( 3.3\right) $
and $\left( 3.10\right) ,$ we have 
\begin{eqnarray}
I_{1} &\leq &C_{1}r^{1-\varepsilon /\mu }C_{0}\left( L_{0}r^{\varepsilon
}\right) ^{\left| \gamma \right| }S\left( k,\gamma \right) L^{k}  \notag \\
\ &\leq &2^{\frac{\sigma \mu }{\omega }+1}C_{1}C_{0}\left(
L_{0}r^{\varepsilon }\right) ^{\left| \gamma \right| }S\left( k+1,\gamma
\right) L^{k}  \label{3.12}
\end{eqnarray}
Using $\left( 3.2\right) $ and $\left( 3.10\right) ,$ we obtain 
\begin{eqnarray*}
I_{2} &\leq &\sum\limits_{\beta <\gamma }\left( _{\beta }^{\gamma }\right)
M^{\left| \gamma -\beta \right| +1}\left[ \Gamma (<\gamma -\beta ,q>+1)%
\right] ^{\sigma \mu }rC_{0}\left( L_{0}r^{\varepsilon }\right) ^{\left|
\beta \right| }S\left( k,\beta \right) L^{k} \\
\ &\leq &\sum\limits_{\beta <\gamma }\left( _{\beta }^{\gamma }\right) \left[
\Gamma (<\gamma -\beta ,q>+1)\right] ^{\sigma \mu }M^{\left| \gamma -\beta
\right| +1}r^{\varepsilon }C_{0}\left( L_{0}r^{\varepsilon }\right) ^{\left|
\beta \right| }2^{\frac{\sigma \mu }{\omega }+1}S\left( k+1,\beta \right)
L^{k}
\end{eqnarray*}
In the other hand, $\left( 3.5\right) $ and $\left( 3.6\right) ,$ give 
\begin{equation}
\left( _{\beta }^{\gamma }\right) \left[ \Gamma (<\gamma -\beta ,q>+1)\right]
^{\sigma \mu }S\left( k,\beta \right) \leq C_{2}^{\sigma \mu <\gamma -\beta
,q>}S\left( k,\gamma \right)  \label{3.13}
\end{equation}
Thus we obtain 
\begin{eqnarray*}
I_{2} &\leq &\sum\limits_{\beta <\gamma }C_{2}^{\sigma \mu \left| \gamma
-\beta \right| }M^{\left| \gamma -\beta \right| +1}r^{\varepsilon
}C_{0}\left( L_{0}r^{\varepsilon }\right) ^{\left| \beta \right| }2^{\frac{%
\sigma \mu }{\omega }+1}S\left( k+1,\gamma \right) L^{k} \\
\ &\leq &\sum\limits_{\beta <\gamma }\left( \frac{MC_{2}^{\sigma \mu }}{%
L_{0}r^{\varepsilon }}\right) ^{\left| \gamma -\beta \right| }r^{\varepsilon
}2^{\frac{\sigma \mu }{\omega }+1}MC_{0}\left( L_{0}r^{\varepsilon }\right)
^{\left| \gamma \right| }S\left( k+1,\gamma \right) L^{k} \\
\ &\leq &\sum\limits_{0<\beta \leq \gamma }\left( \frac{MC_{2}^{\sigma \mu }%
}{L_{0}r^{\varepsilon }}\right) ^{\left| \beta \right| }r^{\varepsilon }2^{%
\frac{\sigma \mu }{\omega }+1}MC_{0}\left( L_{0}r^{\varepsilon }\right)
^{\left| \gamma \right| }S\left( k+1,\gamma \right) L^{k} \\
&\leq &\frac{nMC_{2}^{\sigma \mu }}{L_{0}r^{\varepsilon }}\sum\limits_{\beta
\geq 0}\left( \frac{MC_{2}^{\sigma \mu }}{L_{0}r^{\varepsilon }}\right)
^{\left| \beta \right| }r^{\varepsilon }2^{\frac{\sigma \mu }{\omega }%
+1}MC_{0}\left( L_{0}r^{\varepsilon }\right) ^{\left| \gamma \right|
}S\left( k+1,\gamma \right) L^{k}
\end{eqnarray*}
Set $C_{3}=\sum\limits_{\alpha \geq 0}\left( \frac{1}{2}\right) ^{\left|
\alpha \right| },$ since $L_{0}\geq 2MC_{2}^{\sigma \mu }$ and $r\geq 1,$
then 
\begin{equation}
I_{2}\leq \frac{nMC_{2}^{\sigma \mu }}{L_{0}}C_{3}2^{\frac{\sigma \mu }{%
\omega }+1}MC_{0}\left( L_{0}r^{\varepsilon }\right) ^{\left| \gamma \right|
}S\left( k+1,\gamma \right) L^{k}  \label{3.14}
\end{equation}
Finally in view of $\left( 3.2\right) $%
\begin{eqnarray*}
I_{3} &\leq &\sum\limits_{0<<\alpha ,q>\leq 1}\sum\limits_{\beta \leq \gamma
}\left( _{\beta }^{\gamma }\right) [\Gamma (<\gamma -\beta ,q>+1)]^{\sigma
\mu }M^{\left| \gamma -\beta \right| +1}r^{1-<\alpha ,q>} \\
&&\times C_{0}\left( L_{0}r^{\varepsilon }\right) ^{\left| \beta +\alpha
\right| }S\left( k,\beta +\alpha \right) L^{k}
\end{eqnarray*}
For any $\alpha \in Z_{+}^{n},$ $0<<\alpha ,q>\leq 1,$ we have 
\begin{equation*}
r^{1-<\alpha ,q>+\varepsilon <\alpha ,q>}\ \leq \lambda ^{1-<\alpha ,q>}\tau
^{<\alpha ,q>}\quad ,
\end{equation*}
which gives, with $\left( 3.11\right) $ and $\left( 3.13\right) ,$%
\begin{eqnarray*}
I_{3}\ &\leq &\sum\limits_{0<<\alpha ,q>\leq 1}\sum\limits_{\beta \leq
\gamma }\left( _{\beta }^{\gamma }\right) \left[ \Gamma (<\gamma -\beta
,q>+1)\right] ^{\sigma \mu }M^{\left| \gamma -\beta \right| +1} \\
&&\times C_{0}\left( L_{0}r^{\varepsilon }\right) ^{\left| \beta \right|
}L_{0}^{\left| \alpha \right| }2^{\frac{\sigma \mu }{\omega }+1}S\left(
k+1,\beta \right) L^{k} \\
\ &\leq &\sum\limits_{0<<\alpha ,q>\leq 1}\sum\limits_{\beta \leq \gamma
}\left( \frac{MC_{2}^{\sigma \mu }}{L_{0}r^{\varepsilon }}\right) ^{\left|
\gamma \right| }2^{\frac{\sigma \mu }{\omega }+1}MC_{0}L_{0}^{\left| \alpha
\right| }\left( L_{0}r^{\varepsilon }\right) ^{\left| \gamma \right|
}S\left( k+1,\gamma \right) L^{k}
\end{eqnarray*}
Put $C_{4}=\sum\limits_{0<<\alpha ,q>\leq 1}L_{0}^{\left| \alpha \right| },$
then we obtain 
\begin{equation}
I_{3}\leq 2^{\frac{\sigma \mu }{\omega }+1}MC_{4}C_{3}C_{0}\left(
L_{0}r^{\varepsilon }\right) ^{\left| \gamma \right| }S\left( k+1,\gamma
\right) L^{k}  \label{3.15}
\end{equation}
If we choose 
\begin{equation*}
L\geq 2^{\frac{\sigma \mu }{\omega }+1}\left( C_{1}+\frac{%
nM^{2}C_{2}^{\sigma \mu }}{L_{0}}C_{3}+MC_{3}C_{4}\right) \quad ,
\end{equation*}
we get, from $\left( 3.12\right) ,\left( 3.14\right) $ and $\left(
3.15\right) ,$%
\begin{equation*}
I_{1}+I_{2}+I_{3}\leq C_{0}\left( L_{0}r^{\varepsilon }\right) ^{\left|
\gamma \right| }S\left( k+1,\gamma \right) L^{k+1}\quad ,
\end{equation*}
which means that $\left( 3.9\right) $ holds at the order $k+1.$

\textbf{End of the proof of theorem 3.} Applying the last lemma for $\gamma
=0,$ we find 
\begin{eqnarray}
\left| A_{i_{k}...i_{0}}\left( x,r\right) \right| &\leq &C_{0}L^{k}\left(
r^{\left( 1-\varepsilon /\mu \right) k}+\left[ \Gamma \left( k+1\right) %
\right] ^{\sigma \mu }r^{k\varepsilon \left( 2-1/\mu \right) }\right)  \notag
\\
&\leq &C_{0}^{\prime }L^{k}\left( r^{\left( 1-\varepsilon /\mu \right)
k}+\left( k!\right) ^{\sigma \mu }r^{k\varepsilon \left( 2-1/\mu \right)
}\right)  \label{3.16}
\end{eqnarray}
In the other hand, we know that 
\begin{equation*}
\forall \lambda \in \mathbb{R}_{+},\forall k\in \mathbb{Z}_{+},\forall
s>0,\quad \frac{\left( \frac{\lambda ^{1/s\mu }}{2s\mu }\right) ^{k}}{k!}%
\leq \exp \left( \frac{\lambda ^{1/s\mu }}{2s\mu }\right)
\end{equation*}
So 
\begin{equation*}
\left\{ 
\begin{array}{c}
\lambda ^{k}\leq \left( 2s\mu \right) ^{k\mu s}\left( k!\right) ^{s\mu }\exp
\left( \frac{\lambda ^{1/s\mu }}{2}\right) \qquad \qquad \qquad \\ 
\tau ^{k}\leq \left( 2\left( s-\sigma \right) \mu \right) ^{k\mu \left(
s-\sigma \right) }\left( k!\right) ^{\left( s-\sigma \right) \mu }\exp
\left( \frac{\mu ^{1/\left( s-\sigma \right) \mu }}{2}\right)
\end{array}
\right.
\end{equation*}
Thus we obtain, with $\left( 3.16\right) ,$%
\begin{equation*}
\left| A_{i_{k}...i_{0}}\left( x,r\right) \right| \leq C_{0}^{\prime
}L^{k}\left( 2s\mu \right) ^{k\mu s}\left( k!\right) ^{s\mu }\left[ \exp
\left( \frac{r^{\eta }}{2}\right) +\exp \left( \frac{r^{\eta ^{\prime }}}{2}%
\right) \right] \quad ,
\end{equation*}
where 
\begin{equation*}
\eta ^{\prime }=\frac{\varepsilon \left( 2-1/\mu \right) }{\mu \left(
s-\sigma \right) }\leq \eta =\frac{1-\varepsilon /\mu }{\mu s},
\end{equation*}
since $\varepsilon \leq \dfrac{\mu \left( s-\sigma \right) }{2\mu s-\sigma }%
. $ Therefore 
\begin{equation*}
\left| A_{i_{k}...i_{0}}\left( x,r\right) \right| \leq 2C_{0}^{\prime
}L^{\prime k}\left( k!\right) ^{s\mu }\exp \left( \frac{r^{\eta }}{2}\right)
\end{equation*}
The last estimate gives 
\begin{eqnarray*}
\left| P_{i_{k}...}P_{i_{0}}u\left( x\right) \right| &\leq &2C_{0}^{\prime
}L^{\prime k}\left( k!\right) ^{s\mu }\int_{1}^{+\infty }\exp \left( -\frac{%
r^{\eta }}{2}\right) dr \\
\ &\leq &C^{k+1}\left( k\right) !^{s\mu }\quad ,
\end{eqnarray*}
which means that $u\in G^{s}\left( \Omega ,\left( P_{j}\right)
_{j=1}^{N}\right) $
\end{proof}

\section{Consequences}

The theorems of this work give and unify the results of Bolley-Camus \cite
{B-C} and M\'{e}tivier \cite{M} in the homogenous case, and the results of
Bouzar-Cha\"{i}li \cite{B-Ch1} in the quasi-homogenous case.

\begin{corollary}
Let $\Omega $ be an open subset of $\mathbb{R}^{n}$ and $\sigma >s\geq 1,$
and let $\left( P_{j}\right) _{j=1}^{N}$ be a system of linear differential
operators with coefficients in $G^{q,\sigma }\left( \Omega \right) ,$ then 
\begin{equation*}
\left( P_{j}\right) _{j=1}^{N}\text{ is }q-\text{\textit{quasi-elliptic in} }%
\Omega \Longleftrightarrow G^{s}\left( \Omega ,\left( P_{j}\right)
_{j=1}^{N}\right) \subset G^{q,s}\left( \Omega \right)
\end{equation*}
\end{corollary}

Theorem 2 is also a generalization of the main result of \cite{Z}.

\begin{corollary}
Let $\Omega $\ be an open subset of $\mathbb{R}^{n},s\geq 1$ and let $P$ be
a linear differential operator with coefficients in $G^{\theta ,s}\left(
\Omega \right) ,$ then 
\begin{equation*}
P\text{ is multi-quasi-elliptic in }\Omega \Rightarrow G^{s}\left( \Omega
,P\right) \subset G^{\mathcal{F},s}\left( \Omega \right)
\end{equation*}
\end{corollary}

\begin{remark}
In \cite{Z} the necessity has not been given even for the scalar case.
\end{remark}

\end{document}